\documentclass{amsart}
\usepackage{amsfonts, amsopn, amssymb, mathrsfs}

\newcommand{\continuum}{\ensuremath{\mathfrak{c}}}
\newcommand{\Ck}[1]{\ensuremath{\mathscr{C}({#1})}}
\newcommand{\Czerok}[1]{\ensuremath{\mathscr{C}_0({#1})}}
\newcommand{\dual}[1]{\ensuremath{{#1}^*}}
\newcommand{\ind}[1]{\ensuremath{\mbox{\boldmath{$1$}}_{#1}}}
\newcommand{\mapping}[3]{\ensuremath{{#1}:{#2}\longrightarrow{#3}}}
\newcommand{\norm}[1]{\ensuremath{||{#1}||}}
\newcommand{\normdot}{\ensuremath{||\cdot||}}
\newcommand{\pnorm}[2]{\ensuremath{||{#1}||_{#2}}}
\newcommand{\pnormdot}[1]{\ensuremath{||\cdot||_{#1}}}
\newcommand{\rat}{\mathbb{Q}}
\newcommand{\real}{\mathbb{R}}
\newcommand{\setcomp}[2]{\ensuremath{\{{#1}\;|\;\,{#2}\}}}
\newcommand{\sph}[1]{\ensuremath{S_{#1}}}
\newcommand{\wone}{\ensuremath{\omega_1}}

\DeclareMathOperator{\supp}{supp} \DeclareMathOperator{\ran}{ran}

\newtheorem{thm}{Theorem}
\newtheorem{prop}{Proposition}
\newtheorem{cor}{Corollary}
\newtheorem{lem}{Lemma}
\newtheorem{prob}{Problem}

\theoremstyle{definition}
\newtheorem{defn}[thm]{Definition}

\newtheorem{example}[thm]{Example}

\title{Trees, linear orders and G\^{a}teaux smooth norms}
\author{Richard J. Smith}
\address{Queens' College, Cambridge, CB3 9ET, United Kingdom}
\email{rjs209@cam.ac.uk}

\begin{document}

\begin{abstract}
We introduce a linearly ordered set $Z$ and use it to prove a
necessity condition for the existence of a G\^{a}teaux smooth norm
on $\Czerok{\Upsilon}$, where $\Upsilon$ is a tree. This criterion
is directly analogous to the corresponding equivalent condition for
Fr\'{e}chet smooth norms. In addition, we prove that if
$\Czerok{\Upsilon}$ admits a G\^{a}teaux smooth lattice norm then it
also admits a lattice norm with strictly convex dual norm.
\end{abstract}
\maketitle

\section{Introduction and Preliminaries}

Among the most well-established geometrical properties of norms are
smoothness and strict convexity. A norm $\normdot$ on a Banach space
$X$ is called \textit{G\^{a}teaux smooth}, or just
\textit{G\^{a}teaux}, if, given any $x \in X\backslash\{0\}$, there
exists a functional in $\dual{X}$, denoted by $\norm{x}^\prime$,
such that
\[
\lim_{\lambda \rightarrow 0} \frac{\norm{x + \lambda h} -
\norm{x}}{\lambda} \;=\; \norm{x}^\prime(h)
\]
for all $h \in X$. In addition, if the limit above is uniform for
$h$ in the unit sphere $\sph{X}$, then $\normdot$ is called
\textit{Fr\'{e}chet smooth}, or simply \textit{Fr\'{e}chet}.

Turning now to properties of strict convexity, we say that
$\normdot$ is \textit{strictly convex} if, given $x,y \in X$
satisfying $\norm{x} = \frac{1}{2}\norm{x + y} = \norm{y}$, we have
$x = y$. Of the many stronger cousins of strictly convex norms, we
mention one. The norm $\normdot$ is \textit{locally uniformly
rotund}, or \textit{LUR}, if, given a point $x \in \sph{X}$ and a
sequence $(x_n) \subseteq \sph{X}$ satisfying $\norm{x + x_n}
\rightarrow 2$, we have $\norm{x - x_n} \rightarrow 0$.

Renorming theory is a branch of functional analysis that seeks to
determine the extent to which a given Banach space can be endowed
with equivalent norms sporting certain geometrical properties, such
as the ones above. In this paper, a norm on a given Banach space is
always assumed to be equivalent to the canonical norm. We refer the
reader to \cite{dgz:93} for a comprehensive account of this field up
to 1993, together with the more recent surveys \cite{godefroy:01}
and \cite{zizler:03}.

In recent years, trees have assumed an important role in the field,
both as a source of counterexamples to existing questions and as a
vehicle for exploring new avenues of research; see, for example
\cite{haydon:90}, \cite{haydon:95} and \cite{haydon:99}. We say that
a partially ordered set $(\Upsilon,\preccurlyeq)$ is a \textit{tree}
if, given arbitrary $t \in \Upsilon$, the set of predecessors
$\setcomp{s \in \Upsilon}{s \preccurlyeq t}$, denoted by the
\textit{interval} $(0,t]$, is well-ordered. The set of
\textit{immediate successors} of $t \in \Upsilon$ is denoted by
$t^+$. In this way, trees are a natural generalisation of ordinal
numbers. As well as $(0,t]$, we define the interval $(s,t] =
(0,t]\backslash (0,s]$ for $s \preccurlyeq t$, the \textit{wedge}
$[t,\infty) = \setcomp{u \in \Upsilon}{t \preccurlyeq u}$ and
finally $(t,\infty) = [t,\infty)\backslash\{t\}$. We remark that the
symbols $0$ and $\infty$ are, in this context, convenient notational
devices and not themselves elements of $\Upsilon$.

The scattered locally compact \textit{interval topology} on
$\Upsilon$ is the coarsest topology for which all intervals $(0,t]$
are both open and closed. This topology agrees with the standard
interval topology of any ordinal $\Omega$, if we consider $\Omega$
as a tree. To ensure that this topology is also Hausdorff, we
restrict our attention to trees $\Upsilon$ with the property that
every non-empty, linearly ordered set in $\Upsilon$ has at most one
minimal upper bound. With this topology in mind, we consider the
Banach space $\Czerok{\Upsilon}$ of continuous real-valued functions
vanishing at infinity, and the dual space of measures. We remark
that as $\Upsilon$ is scattered, the weak topology and the topology
of pointwise convergence agree on norm-bounded subsets of
$\Czerok{\Upsilon}$.

Trees and linearly ordered sets enjoy close ties. For a
comprehensive review of these relationships, we refer the reader to
\cite{tod:84}. Given partial orders $P$ and $Q$, we say that the map
$\mapping{\rho}{P}{Q}$ is called \textit{increasing} (respectively
\textit{strictly increasing}) if $\rho(s) \preccurlyeq \rho(t)$
(respectively $\rho(s) \prec \rho(t)$) whenever $s \prec t$.
\textit{Decreasing} and \textit{strictly decreasing} functions are
defined analogously. If there exists a strictly increasing map from
$P$ to a linear order $Q$, we say that $P$ is
$Q$\textit{-embeddable}, or $P \preccurlyeq Q$. Evidently, in this
context, $\preccurlyeq$ is a transitive relation on the class of
partial orders. In much of what follows, $P$ will be a tree and $Q$
a linear order. It is well known that $\Upsilon \preccurlyeq \rat$
if and only if $\Upsilon$ is \textit{special}, which means that
$\Upsilon$ can be written as a countable union of antichains (cf.\
\cite[Theorem 9.1]{tod:84}). Special trees tend to have very good
properties; for example, the following result can be found in
\cite{smith:05b}.

\begin{thm}
\label{speciallur} Given a tree $\Upsilon$, the space
$\Czerok{\Upsilon}$ admits a norm with LUR dual norm if and only if
$\Upsilon$ is special.
\end{thm}

We introduce a couple of combinatorial ideas used extensively in
\cite{haydon:99}.

\begin{defn}
\label{badpoints} Given an increasing function
$\mapping{\rho}{\Upsilon}{\real}$, we say that $t \in \Upsilon$ is a
\textit{bad point for} $\rho$ if there exists a sequence of distinct
points $(u_n) \subseteq t^+$, such that $\rho(u_n) \rightarrow
\rho(t)$.
\end{defn}

Bad points are so named because their presence often indicates that
the given $\Czerok{\Upsilon}$ space has negative renorming
properties. An analogue of the next simple result appears at the
beginning of Section \ref{examples}.

\begin{prop}[(Haydon)]
\label{ratbadpoints} The tree $\Upsilon$ is special if and only if
$\Upsilon \preccurlyeq \real$ and there exists an increasing map
$\mapping{\rho}{\Upsilon}{\real}$ that has no bad points.
\end{prop}

We move on to the second combinatorial property taken from
\cite{haydon:99}.

\begin{defn}
\label{everbranching} A subset $E$ of a tree is said to be
\textit{ever-branching} if each element of $E$ has a pair of strict
successors in $E$ that are incomparable in the tree order.
\end{defn}

It is easy to see that within every ever-branching subset can be
found a \textit{dyadic tree of height} $\omega$; that is, a tree
with a single minimal element, no limit elements, and with the
property that each element has exactly two immediate successors.

Many types of norm on $\Czerok{\Upsilon}$ can be characterised in
terms of increasing real-valued functions on $\Upsilon$, with
further combinatorial properties that can be expressed in terms of
bad points and ever-branching subsets. Of particular interest to us
is the following result.

\begin{thm}[(Haydon \cite{haydon:99})]
\label{f} Given a tree $\Upsilon$, the space $\Czerok{\Upsilon}$
admits a Fr\'{e}chet norm if and only if there exists an increasing
function $\mapping{\rho}{\Upsilon}{\real}$ that has no bad points
and is not constant on any ever-branching subset.
\end{thm}

In order to exhibit a tree that does not satisfy the statement of
Theorem \ref{f}, we introduce a fundamental construction, due to
Kurepa. Given a linear order $\Sigma$, we define the Hausdorff tree
\[
\sigma \Sigma \;=\; \setcomp{A \subseteq \Sigma}{A \mbox{ is
well-ordered}}.
\]
We remark that some authors demand the additional requirement that
elements of $\sigma \Sigma$ are bounded above. One of the reasons
why Kurepa's construction is so important in the theory of trees is
summed up by the following theorem.

\begin{thm}[(Kurepa \cite{kurepa:56})]
\label{sigmanoembed} If $\Sigma$ is a linear order then $\sigma
\Sigma \not \preccurlyeq \Sigma$.
\end{thm}

From Theorem \ref{sigmanoembed}, $\sigma \rat$ is not special. On
the other hand, if we take an enumeration $(q_n)$ of the rationals
and consider the map $A \mapsto \sum_{q_n \in A} 2^{-n}$, we see
that $\sigma \rat \preccurlyeq \real$. It follows that, by
Proposition \ref{ratbadpoints}, every increasing, real-valued
function defined on $\sigma \rat$ has a bad point.

\begin{cor}[(Haydon)]
\label{haydonfthm} The space $\Czerok{\sigma \rat}$ admits no
Fr\'{e}chet norm.
\end{cor}

While many types of norm are accounted for in \cite{haydon:99},
equivalent conditions for the existence of norms on
$\Czerok{\Upsilon}$ with strictly convex dual, or G\^{a}teaux norms,
cannot be adequately expressed in terms of increasing real-valued
functions. In all that follows, $\wone$ denotes the first
uncountable ordinal. The following linearly ordered set is
introduced in \cite{smith:05b}.

\begin{defn}
\label{ordery} Let $Y$ be the set of all strictly increasing,
continuous, transfinite sequences $x = (x_\xi)_{\xi \leq \beta}$ of
real numbers, where $0 \leq \beta < \wone$. Order $Y$ by declaring
that $x < y$ if and only if either $y$ strictly extends $x$, or if
there is some ordinal $\alpha$ such that $x_\xi = y_\xi$ for $\xi <
\alpha$ and $y_\alpha < x_\alpha$.
\end{defn}

Observe that $Y$ is not ordered in the usual lexicographic way.
Compared to the real line, $Y$ is large.

\begin{prop}[(Smith {\cite{smith:05b}})]
\label{ybetaembed} If $\beta < \wone$ then $Y^\beta \preccurlyeq Y$,
where $Y^\beta$ is ordered lexicographically.
\end{prop}

As $\real \preccurlyeq Y$, we see that $\real^\beta \preccurlyeq Y$
for all $\beta < \wone$. On the other hand, it can be shown that $Y$
contains no well-ordered or conversely well-ordered subsets. The
next theorem is the main result of \cite{smith:05b}.

\begin{thm}[(Smith {\cite{smith:05b}})]
\label{dualrotundthm} Given a tree $\Upsilon$, the Banach space
$\Czerok{\Upsilon}$ admits a norm with strictly convex dual norm if
and only if $\Upsilon \preccurlyeq Y$.
\end{thm}

Theorem \ref{dualrotundthm} is a direct analogue of Theorem
\ref{speciallur}. In \cite{smith:05b}, it is shown that the spaces
$\Czerok{\sigma (\real^\beta)}$, where $\real^\beta$ is ordered
lexicographically, admit norms with strictly convex duals provided
$\beta < \wone$. On the other hand, by Theorem \ref{sigmanoembed},
$\Czerok{\sigma Y}$ does not admit such a norm.

The order $Y$ can also be used to give an improved sufficient
condition for the existence of G\^{a}teaux norms in the context of
trees.

\begin{thm}[(Smith {\cite{smith:05c}})]
\label{gsuff} If there exists an increasing function
$\mapping{\rho}{\Upsilon}{Y}$ that is not constant on any
ever-branching subset then $\Czerok{\Upsilon}$ admits a G\^{a}teaux
norm.
\end{thm}

We end our review of the existing literature by presenting what was
hitherto the best known necessary condition for G\^{a}teaux norms in
this context. Given a tree $\Upsilon$, the \textit{forcing topology}
on $\Upsilon$ takes as its basis the set of all wedges $[t,\infty)$,
$t \in \Upsilon$. A subset $B \subseteq \Upsilon$ is called
\textit{Baire} if it is a Baire space with respect to the induced
forcing topology; that is, any countable intersection of relatively
dense, open subsets of $B$ is again dense. When referring to the
Baire property, we will only consider subsets that are
\textit{perfect} with respect to the forcing topology; in other
words those without isolated points or, equivalently, maximal
elements. Arguably the simplest example of such an object is the
ordinal $\wone$, though more interesting ones that have no
uncountable linearly ordered subsets can be found in \cite[Lemma
9.12]{tod:84} (cf.\ \cite{haydon:95}).

Theorems \ref{f} and \ref{gsuff} applied to a constant function on
$\wone$ demonstrate that, by itself, the Baire property cannot
destroy G\^{a}teaux renormability. Instead, we have the following
result.

\begin{thm}[(Haydon {\cite{haydon:95}})]
\label{haydongnec} If $\Czerok{\Upsilon}$ admits a G\^{a}teaux norm
then $\Upsilon$ contains no ever-branching Baire subsets.
\end{thm}

We turn now to the results of this paper. In order to properly
express our necessary condition for G\^{a}teaux renormability, we
must introduce a second linearly ordered set.

\begin{defn}
\label{orderz} Let $Z$ be the set of all increasing, continuous
sequences $x = (x_\xi)_{\xi \leq \beta}$ of real numbers, where $0
\leq \beta < \wone$, and such that $x$ is strictly increasing on
$[0,\beta)$. The order of $Z$ follows that of $Y$; $x < y$ if and
only if either $y$ strictly extends $x$, or if there is some ordinal
$\alpha$ such that $x_\xi = y_\xi$ for $\xi < \alpha$ and $y_\alpha
< x_\alpha$.
\end{defn}

The elements of $Z$ that are not in $Y$ are exactly those of the
form $x = (x_\xi)_{\xi \leq \beta+1}$, where $(x_\xi)_{\xi \leq
\beta} \in Y$ and $x_\beta = x_{\beta+1}$. This order is a partial
Dedekind completion of $Y$. We also need a natural definition of bad
points with respect to $Z$.

\begin{defn}
\label{zbadpoint} Given an increasing function
$\mapping{\rho}{\Upsilon}{Z}$, we say that $t \in \Upsilon$ is
\textit{$Z$-bad} for $\rho$ if there exists a sequence of distinct
points $(u_n) \subseteq t^+$ such that $\rho(u_n) \rightarrow
\rho(t)$ in the order topology of $Z$.
\end{defn}

Using $Z$-bad points, we obtain a direct analogy to the necessity
part of Theorem \ref{f}; the following is the main result of this
paper.

\begin{thm}
\label{newgnec} If the space $\Czerok{\Upsilon}$ admits a
G\^{a}teaux norm, then there exists an increasing function
$\mapping{\rho}{\Upsilon}{Z}$ that has no $Z$-bad points and is not
constant on any ever-branching subset.
\end{thm}

In some sense, $Y$ is to $\rat$ what $Z$ is to $\real$, and these
relationships correspond well to those of Theorems
\ref{dualrotundthm}, \ref{speciallur}, \ref{newgnec} and \ref{f}
respectively.

The following corollary of Theorem \ref{newgnec} generalises a
result from \cite{fhz:97}, which states that $\Czerok{[0,\wone)}$
does not admit any G\^{a}teaux lattice norm.

\begin{cor}
\label{latticeggivesr} If $\Czerok{\Upsilon}$ admits a G\^{a}teaux
lattice norm then $\Upsilon \preccurlyeq Y$ and, consequently,
$\Czerok{\Upsilon}$ admits a lattice norm with strictly convex dual.
\end{cor}

We end Section \ref{necessitycondition} by proving the next
proposition, which shows that Theorem \ref{haydongnec} is a
corollary of Theorem \ref{newgnec}.

\begin{prop}
\label{haydoncor} If $\mapping{\rho}{\Upsilon}{Z}$ is an increasing
function that is not constant on any ever-branching subset, then
$\Upsilon$ does not admit any ever-branching Baire subsets.
\end{prop}

The final section, devoted to examples, begins with a proof that
Theorem \ref{haydongnec} is strictly implied by Theorem
\ref{newgnec}.

\begin{prop}
\label{sigmaYnoG} The tree $\sigma Y$ is $Z$-embeddable, but every
increasing function \mapping{\rho}{\Upsilon}{Z} has a $Z$-bad point.
In particular, $\Czerok{\sigma Y}$ does not admit a G\^{a}teaux
norm.
\end{prop}

Proposition \ref{sigmaYnoG} is analogous to Corollary
\ref{haydonfthm}. Section \ref{examples} ends with Example
\ref{Ggap}, which shows that there is a gap between the conditions
of Theorems \ref{gsuff} and \ref{newgnec}. This, together with the
analogies presented above and the author's bias, prompts the
following problem.

\begin{prob}
\label{gateauxconj} If there exists an increasing function
$\mapping{\rho}{\Upsilon}{Z}$ that has no $Z$-bad points and is not
constant on any ever-branching subset, does $\Czerok{\Upsilon}$
admit a G\^{a}teaux norm?
\end{prob}

Recently, the author gave a purely topological formulation of
Theorem \ref{dualrotundthm}. Given a tree $\Upsilon$, the space
$\Czerok{\Upsilon}$ admits a norm with strictly convex dual norm if
and only if $\Upsilon$ is a so-called \textit{Gruenhage space}, with
respect to its interval topology \cite{smith:07}.

\begin{prob}
\label{internal} Is there an internal characterisation of trees
$\Upsilon$, with the property that $\Czerok{\Upsilon}$ admits a
G\^{a}teaux norm?
\end{prob}

Problem \ref{internal} may be restated in terms of Fr\'{e}chet
norms, Kadec norms and others. This section closes with further
problem, motivated by Corollary \ref{latticeggivesr}.

\begin{prob}
\label{latticegdualrgeneral} If $L$ is locally compact and
$\Czerok{L}$ admits a G\^{a}teaux lattice norm, does $\Czerok{L}$
admit a norm with strictly convex dual? Is this statement also true
with respect to a general Banach lattice?
\end{prob}

\section{Necessity conditions for G\^{a}teaux renormability}
\label{necessitycondition}

To help familiarise the reader with $Z$ and $Z$-bad points, we begin
by briefly describing some forms of sequential convergence in $Z$.
First observe that if $x \in Y$, $y \in Z$ and $y > x$ is
sufficiently close to $x$ in the order topology of $Z$, then $y$
must be a strict extension of $x$. On the other hand, if $x \in
Z\backslash Y$ then $x$ has no strict extensions in $Z$. The proof
of the next lemma is a simple exercise in elementary analysis and is
omitted.

\begin{lem}
\label{convergenceinZ} Let $x \in Z$ and suppose $(z^n) \subseteq Z$
is a sequence satisfying $x < z^n$. We have the following rules for
the convergence of $(z^n)$ to $x$:
\begin{enumerate}
\item[1.] if $x = (x_\xi)_{\xi \leq \beta} \in Y$ then
$z^n \rightarrow x$ if and only if $z^n$ strictly extends $x$ for
large enough $n$, and $z^n_{\beta + 1} \rightarrow \infty$.
\end{enumerate}
If $x = (x_\xi)_{\xi \leq \beta+1} \in Z\backslash Y$ then since $x$
has no strict extensions, there exists $\alpha_n \leq \beta$ such
that $z^n_\xi = x_\xi$ for $\xi < \alpha_n$ and $z^n_{\alpha_n} <
x_{\alpha_n}$. In this case, we have:
\begin{enumerate}
\item[2.] if $\beta = 0$ or $\beta = \alpha + 1$ for some $\alpha$, then $z^n \rightarrow x$ if and only
if $\alpha_n = \beta$ for large enough $n$, and $z^n_\beta
\rightarrow x_\beta$;
\item[3.] if $\beta$ is a limit ordinal, then $z^n \rightarrow x$ if and only
if $\alpha_n \rightarrow \beta$.
\end{enumerate}
\end{lem}

We present a simple application of Lemma \ref{convergenceinZ}. If
$\mapping{\pi}{\Upsilon}{Y}$ is a strictly increasing map then it
could have $Z$-bad points. However, if we fix an order isomorphism
$\mapping{\theta}{\real}{(0,1)}$ and define, for $x = (x_\xi)_{\xi
\leq \beta} \in Y$, $\Theta(x)_\xi = \theta(x_\xi)$ whenever $\xi
\leq \beta$, then by Lemma \ref{convergenceinZ} part (1), the
strictly increasing $Y$-valued map $\Theta \circ \pi$ has no $Z$-bad
points. Thus, some $Z$-bad points are easily removed by making
simple adjustments. More details of how $Z$ operates can be found in
Section \ref{examples}.

Now, for the rest of this section, we fix a norm $\normdot$ on
$\Czerok{\Upsilon}$. We continue by introducing a concept that
features in both \cite{haydon:95} and \cite{haydon:99}. Given $t \in
\Upsilon$, let $C_t$ be the set of all $f \in \Czerok{\Upsilon}$
such that $f$ vanishes outside $(0,t]$ and increasing on $(0,t]$.

\begin{defn}
\label{mu-function} If $f \in C_t$ and $\delta \geq 0$, the
increasing function $\mu(f,\delta,\cdot)$ is defined on the wedge
$[t,\infty)$ by
\[
\mu(f,\delta,\cdot) \;=\; \inf \setcomp{\norm{f + (f(t) + \delta)
\ind{(t,u]} + \varphi}}{\varphi \in \Czerok{\Upsilon} \mbox{ and }
\supp \varphi \subseteq (u,\infty)}
\]
where $\ind{A}$ denotes the indicator function of the set $A$ and
$\supp \varphi$ is the support of $\varphi$. We also define the
abbreviation $\mu(f,\cdot)$ by $\mu(f,u) = \mu(f,0,u)$ and the
associated function $\mu$, given by $\mu(t) = \inf
\setcomp{\norm{\ind{(0,t]} + \varphi}}{\varphi \in \Czerok{\Upsilon}
\mbox{ and } \supp \varphi \subseteq (t,\infty)}$.
\end{defn}

Attainment of the infimum in the definition of these so-called
$\mu$-functions has important consequences for the renormability of
$\Czerok{\Upsilon}$, and bad points and ever-branching subsets come
into play. The first consequence of the following lemma is trivial,
and the second and third are immediate generalisations of
\cite[Lemma 3.1]{haydon:99} and \cite[Proposition 3.4]{haydon:99}
respectively.

\begin{lem}[(Haydon \cite{haydon:99})]
\label{infattaining} Suppose $t \in \Upsilon$, $f \in C_t$ and
$\delta \geq 0$. Then:
\begin{enumerate}
\item if $\normdot$ is a lattice norm then $\norm{f + (f(t) + \delta)
\ind{(t,u]}} = \mu(f,\delta,u)$ for all $u \succcurlyeq t$;
\item if $u \succcurlyeq t$ is a bad point for $\mu(f,\delta,\cdot)$ then
$\norm{f + (f(t) + \delta) \ind{(t,u]}} = \mu(f,\delta,u)$;
\item if $\mu(f,\delta,\cdot)$ is constant on some ever-branching
subset $E \subseteq (u, \infty)$, where $u \succcurlyeq t$, then
there exists $\varphi \in \Czerok{\Upsilon}$ with
\[
\supp \varphi \;\subseteq\; \setcomp{v \in (u,\infty)}{v
\preccurlyeq w \mbox{ for some } w \in E}
\]
and $\mu(f,\delta,u) = \norm{f + (f(t) + \delta)(\ind{(t,u]} +
\varphi)}$.
\end{enumerate}
\end{lem}

We continue with an idea from \cite{smith:05b}.

\begin{defn}
\label{plateau} A subset $V \subseteq \Upsilon$ is called a
\textit{plateau} if $V$ has a least element $0_V$ and $V =
\bigcup_{t \in V}[0_V,t]$. A partition $\mathscr{P}$ of $\Upsilon$
consisting solely of plateaux is called a \textit{plateau
partition}.
\end{defn}

Observe that if $V$ is a plateau then $V\backslash\{0_V\}$ is open.
It follows that if we have a plateau partition $\mathscr{P}$ and
define the \textit{set of least elements} $H = \setcomp{0_V}{V \in
\mathscr{P}}$, then $H$ is closed in $\Upsilon$. Of course, $H$ may
be regarded as a tree in its own right, with its own interval
topology. Plateaux are stable under taking arbitrary intersections.

\begin{prop}[(Smith {\cite[Proposition 10]{smith:05b}})]
\label{plateauintersection} Let $\Upsilon$ be a tree and
$\mathfrak{F}$ a family of plateaux of $\Upsilon$ with non-empty
intersection $W$. Then $W$ is a plateau and $0_W = \sup_{V \in
\mathfrak{F}} 0_V$.
\end{prop}

The connection between increasing functions and plateaux is given by
the next proposition.

\begin{prop}[(Smith {\cite[Proposition 9]{smith:05b}})]
\label{plateaupartition} Let $\mapping{\rho}{\Upsilon}{\Sigma}$ be
an increasing function into a linear order $\Sigma$. Then the
equivalence relation $\sim$, given by $s \sim t$ if and only if
there exists $r \preccurlyeq s,t$ such that $\rho(s) = \rho(r) =
\rho(t)$, defines the plateau partition of $\Upsilon$, with respect
to $\rho$. Moreover, the restriction of $\rho$ to the set of least
elements $H = \setcomp{0_V}{V \in \mathscr{P}}$ is strictly
increasing.
\end{prop}

Proposition \ref{plateaupartition} applies equally well to
decreasing functions. As the $\mu$-functions from Definition
\ref{mu-function} are increasing on their respective domains, they
may be analysed using plateaux. Elements of the following technical
lemma appear implicitly in the proof of \cite[Theorem
8.1]{haydon:99}.

\begin{lem}
\label{lambda} Let $\normdot$ be G\^{a}teaux smooth and suppose that
$\varepsilon\pnormdot{\infty} \leq \normdot \leq \pnormdot{\infty}$
for some $\varepsilon \in (0,1)$. Moreover, suppose $V$ is a
plateau, $f \in C_{0_V}$ and $\mu(f,\cdot)$ is constant on $V$. We
define a function $\lambda$ on $V\backslash\{0_V\}$ by setting
\[
\lambda(t) \;=\; \sup \setcomp{\delta \geq 0}{\mu(f,\delta,t) \leq
\mu(f,0_V) + \textstyle{\frac{1}{2}}\varepsilon\delta}.
\]
We check that $\lambda$ is well-defined and satisfies the following
properties:
\begin{enumerate}
\item $\lambda$ is decreasing on $V\backslash\{0_V\}$;
\item if $\lambda$ takes constant value $\nu$ on the plateau $W \subseteq
V\backslash\{0_V\}$ then $\mu(f,\nu,\cdot)$ takes constant value
$\mu(f,0_V) + \frac{1}{2}\varepsilon\nu$ on $W$;
\item if $\mathscr{P}$ is the plateau partition of $V\backslash\{0_V\}$
with respect to $\lambda$, supplied by Proposition
\ref{plateaupartition}, $W \in \mathscr{P}$, and $f_W \in C_{0_W}$
is defined by
\[
f_W \;=\; f + (f(0_V) + \lambda(0_W))\ind{(0_V,0_W]}
\]
then $\mu(f_W,\cdot)$ takes constant value $\mu(f,0_V) +
\frac{1}{2}\varepsilon \lambda(0_W)$ on $W$;
\item if the infimum in the definition of $\mu(f,t)$ is attained
then $\lambda(t) > 0$.
\end{enumerate}
\end{lem}

\begin{proof}
Fix $t \in V\backslash\{0_V\}$ and, for $\delta \geq 0$, define
$F(\delta) = \mu(f,\delta,t) - \mu(f,0_V) -
\frac{1}{2}\varepsilon\delta$. Observe that $F$ is continuous and
$F(0) = 0$. Moreover, if $\supp \varphi$ is a subset of
$(t,\infty)$, we estimate that $\norm{f + (f(t) +
\delta)\ind{(0_V,t]} + \varphi} \geq \varepsilon\delta - \norm{f +
f(t)\ind{(0_V,t]}}$, whence $F(\delta)$ tends to $\infty$ as
$\delta$ does. As a result, $\lambda(t)$ is well-defined.

Now we can check the properties of $\lambda$. We see that
$\mu(f,\lambda(t),t) = \mu(f,0_V) +
\frac{1}{2}\varepsilon\lambda(t)$ for any $t \in
V\backslash\{0_V\}$. Therefore, if $t \preccurlyeq u$ then, as
$\mu(f,\lambda(u),\cdot)$ is increasing, we have
\[
\mu(f,\lambda(u),t) \;\leq\; \mu(f,\lambda(u),u) \;=\; \mu(f,0_V) +
\textstyle{\frac{1}{2}}\varepsilon\lambda(u)
\]
which shows that $\lambda(t) \geq \lambda(u)$, giving us property
(1).

The second property follows immediately and the third follows from
the second. To prove property (4), we let $g = f + f(t)\ind{(0_V,t]}
+ \varphi$ with $\supp \varphi \subseteq (t,\infty)$, such that
$\norm{g} = \mu(f,t) = \mu(f,0_V)$. Observe that as the infimum
$\mu(f,0_V)$ is attained, we have
\[
\norm{g}^\prime(\ind{(0_V,t]}) \;=\; \lim_{\delta \rightarrow 0_+}
\frac{\norm{g + \delta \ind{(0_V,t]}} - \norm{g}}{\delta} \;\geq\; 0
\]
and similarly for $-\ind{(0_V,t]}$, whence
$\norm{g}^\prime{(\ind{(0_V,t]})} = 0$. Now it is evident that there
exists $\delta > 0$ satisfying
\[
\mu(f,\delta,t) \;\leq\; \norm{g + \delta\ind{(0_V,t]}} \;\leq\;
\norm{g} + \textstyle{\frac{1}{2}}\varepsilon\delta \;=\; \mu(f,0_V)
+ \textstyle{\frac{1}{2}}\varepsilon\delta
\]
which means that $\lambda(t) \geq \delta > 0$.
\end{proof}

While noting property (4) above, we stress that sometimes $\lambda$
does vanish, and it is necessary to analyse what happens in this
case.

\begin{lem}
\label{lambdavanish} Suppose $V$, $f$, $\mu(f,\cdot)$, $\lambda$ and
the partition $\mathscr{P}$ are as in Lemma \ref{lambda}. If
$\lambda(t) = 0$ for some $t \in W \in \mathscr{P}$, then:
\begin{enumerate}
\item $W = [0_W,\infty) \cap V$;
\item $W$ is finitely-branching, in other words, $u^+ \cap W$ is
finite whenever $u \in W$;
\item $W$ contains no ever-branching subsets.
\end{enumerate}
\end{lem}

\begin{proof}
The first property follows because $\lambda \geq 0$ and is
decreasing. To prove property (2), we suppose that $u \in V$ is such
that $u^+ \cap V$ is infinite. Then $u$ is a bad point for
$\mu(f,\cdot)$ as $\mu(f,v) = \mu(f,u)$ for infinitely many $v \in
u^+$. Consequently, the infimum in the definition of $\mu(f,u)$ is
attained by part (2) of Lemma \ref{infattaining}, and it follows
from Lemma \ref{lambda} part (4) that $\lambda(u) > 0$. As a result,
$u \notin W$. For property (3), it is enough to show that if $u \in
V$ and $E$ is an ever-branching subset of $[u,\infty) \cap V$, then
$\lambda(u)
> 0$. Indeed, given such $u$ and $E$, by part (3) of Lemma \ref{infattaining},
the infimum in the definition of $\mu(f,u)$ is attained. Therefore,
by part (4) of Lemma \ref{lambda}, $\lambda(u)
> 0$.
\end{proof}

The proof of Theorem \ref{newgnec} is similar to that of Theorem
\ref{dualrotundthm}, in that it employs monotone real-valued
functions to recursively define a refining sequence of plateaux
partitions of the given tree. This sequence is used to define a
$Z$-valued function or, in the case of Theorem \ref{dualrotundthm}
or Corollary \ref{latticeggivesr}, a $Y$-valued function. We will
see that we must make use of the elements in $Z\backslash Y$
precisely when our $\lambda$-functions from Lemma \ref{lambda}
vanish.

\begin{proof}[of Theorem \ref{newgnec}]
Let $\normdot$ be G\^{a}teaux smooth and suppose that
$\varepsilon\pnormdot{\infty} \leq \normdot \leq \pnormdot{\infty}$
for some $\varepsilon \in (0,1)$. We assemble, for each $\beta <
\wone$, a plateau partition $\mathscr{P}_{\beta}$, and for each $V
\in \mathscr{P}_\beta$, a function $f_{(\beta,V)} \in C_{0_V}$ such
that:
\begin{enumerate}
\item $\mu(f_{(\beta,V)},\cdot)$ takes constant value
$\mu(f_{(\beta,V)},0_V)$ on $V$;
\item $\mu(f_{(\beta,V)},0_V) - 1 \;\leq\;
\frac{1}{2}\varepsilon(\pnorm{f_{(\beta,V)}}{\infty} - 1)$.
\end{enumerate}
Following this, we define a function $\mapping{\pi}{\Upsilon}{Z}$
and prove that it possesses a number of properties. Our final
function $\rho$ will be a modification of $\pi$.

We begin by constructing $\mathscr{P}_0$. Recall the increasing
function $\mu$ from Definition \ref{mu-function}. Let
$\mathscr{P}_0$ be its plateau partition, courtesy of Proposition
\ref{plateaupartition}, and define $f_{(0,V)} = \ind{(0,0_V]}$ for
$V \in \mathscr{P}_0$. It follows that $\mu(f_{(0,V)},\cdot)$ takes
constant value $\mu(f_{(0,V)},0_V) = \mu(0_V)$ on $V$, and that
\[
\mu(f_{(0,V)},0_V) - 1 \;\leq\; \norm{\ind{(0,0_V]}} - 1 \;\leq\; 0
\;=\; \textstyle{\frac{1}{2}}\varepsilon(\pnorm{f_{(0,V)}}{\infty} -
1).
\]

Now suppose $\mathscr{P}_\beta$ and the associated $f_{(\beta,V)}$
have been built. Let $V \in \mathscr{P}_\beta$. If $V = \{0_V\}$
then set $\mathscr{P}_V = \{V\}$ and $f_{(\beta+1,V)} =
f_{(\beta,V)}$. Otherwise, Lemma \ref{lambda}, together with
Proposition \ref{plateaupartition}, furnishes us with the plateau
partition of $V\backslash\{0_V\}$ associated with the
$\lambda$-function. We augment this with the single element
$\{0_V\}$ to give a plateau partition $\mathscr{P}_V$ of $V$. Set
$\mathscr{P}_{\beta+1} = \bigcup\setcomp{\mathscr{P}_V}{V \in
\mathscr{P}_\beta}$. If $W \in \mathscr{P}_V$ then either $W =
\{0_V\}$ or $W \subseteq V\backslash\{0_V\}$. In the former case let
$f_{(\beta+1,W)} = f_{(\beta,V)}$; it is easy to see that
$f_{(\beta+1,W)}$ satisfies conditions (1) and (2) above. In the
latter case, let $f_{(\beta+1,W)} = f_W$, where $f_W$ is as in Lemma
\ref{lambda} part (3). We observe condition (1) is satisfied, again
by Lemma \ref{lambda} part (3). To see that condition (2) holds,
note that
\[
\mu(f_{(\beta+1,W)},0_W) - \mu(f_{(\beta,V)},0_V) \;=\;
\textstyle{\frac{1}{2}}\varepsilon\lambda(0_W) \;=\;
\textstyle{\frac{1}{2}}\varepsilon(\pnorm{f_{(\beta+1,W)}}{\infty} -
\pnorm{f_{(\beta,V)}}{\infty})
\]
and apply the inductive hypothesis.

We move on to the limit case. Suppose that $\beta < \wone$ is a
limit ordinal and that all has been constructed for $\alpha <
\beta$. Given $t \in \Upsilon$, we let $V_\alpha^t \in
\mathscr{P}_\alpha$ be such that $t \in V_\alpha^t$. Set
$\mathscr{P}_\beta = \setcomp{\bigcap_{\alpha < \beta} V_\alpha^t}{t
\in \Upsilon}$. Fix some $V \in \mathscr{P}_\beta$. Let $t = 0_V$,
$V_\alpha = V_\alpha^t$, $t_\alpha = 0_{V_\alpha}$ and $f_\alpha =
f_{(\alpha,V_\alpha)}$. Then $t = \sup_{\alpha < \beta} t_\alpha$ by
Proposition \ref{plateauintersection}. What we would like to do is
define $f_{(\beta,V)} = f \in \Czerok{\Upsilon}$ to be the unique
function supported on $(0,t]$, such that its restriction to
$(0,t_\alpha]$ is $f_\alpha$. This can indeed be done, provided that
$(\pnorm{f_\alpha}{\infty})_{\alpha < \beta}$ is bounded. Observe
that if $g \in C_u$ satisfies condition (2) above then
\[
\varepsilon \pnorm{g}{\infty}-1 \;\leq\; \mu(g,u)-1 \;\leq\;
\textstyle{\frac{1}{2}}\varepsilon(\pnorm{g}{\infty}-1)
\]
giving $\pnorm{g}{\infty} \leq \frac{2}{\varepsilon} - 1$. Therefore
$(\pnorm{f_\alpha}{\infty})_{\alpha < \beta}$ is bounded as
required. Moreover, since each $f_\alpha \in C_{t_\alpha}$, we have
$f \in C_t$. Now set $g_\alpha = f_\alpha +
f_\alpha(t_\alpha)\ind{(t_\alpha,t]}$. Of course, as $f_\alpha$ is
increasing on $(0,t_\alpha]$ and vanishes elsewhere, we have
$\pnorm{g_\alpha}{\infty} = \pnorm{f_\alpha}{\infty}$. Moreover, as
$\mu(f_\alpha,\cdot)$ takes constant value $\mu(f_\alpha,t_\alpha)$
on $V_\alpha$ by inductive hypothesis, and $\mu(g_\alpha,u) =
\mu(f_\alpha,u)$ whenever $u \in V \subseteq V_\alpha$, it follows
that $\mu(g_\alpha,\cdot)$ takes constant value
$\mu(f_\alpha,t_\alpha)$ on $V$. The reader can verify that, as
$(g_\alpha)_{\alpha < \beta}$ converges in norm to $f$,
$(\mu(g_\alpha,\cdot))_{\alpha < \beta}$ converges uniformly to
$\mu(f,\cdot)$ (cf.\ \cite[Lemma 3.6]{haydon:99}). As a result, $f$
satisfies conditions (1) and (2) above. This ends the recursion.

Now we define $\pi$. Given $t \in \Upsilon$, let $V^t_\beta$ be as
above. In addition, we let $\lambda^t_\beta$ be the
$\lambda$-function associated with $V^t_\beta$ and
$f_{(\beta,V^t_\beta)}$, provided $V^t_\beta$ is not a singleton.
Set $\pi(t)_0 = -\mu(t)$. If $\beta > 0$, let $\pi(t)_\beta =
\mu(f_{(\beta,V^t_\beta)},t)$ as long as $0_{V^t_\alpha} \prec t$
for all $\alpha < \beta$ and $\lambda^t_\alpha(t) > 0$ whenever
$\alpha+1 < \beta$. Otherwise, we leave $\pi(t)_\beta$ undefined.

We verify that $\pi(t)$ is an element of $Z$. Observe that if
$\pi(t)_\beta$ is defined, then so is $\pi(t)_\alpha$ whenever
$\alpha < \beta$. If $0 < \alpha < \beta$ then $\pi(t)_0 < 0 <
\pi(t)_\alpha$ and moreover
\begin{eqnarray*}
\pi(t)_{\alpha+1} &=& \mu(f_{(\alpha+1,V^t_{\alpha+1})},t) \\
&=& \mu(f_{(\alpha,V^t_\alpha)},t) +
\textstyle{\frac{1}{2}}\varepsilon\lambda^t_\alpha (0_{V^t_{\alpha+1}}) \\
&=& \pi(t)_\alpha + \textstyle{\frac{1}{2}}\varepsilon
\lambda^t_\alpha(t)
\end{eqnarray*}
whence $\pi(t)_{\alpha+1} \geq \pi(t)_\alpha$. In addition, if
$\alpha + 1 < \beta$ then $\pi(t)_{\alpha+1} > \pi(t)_\alpha$ by our
definition of $\pi$. Now, if $\beta$ is a limit ordinal and
$\pi(t)_\alpha$ is defined for all $\alpha < \beta$, so is
$\pi(t)_\beta$. Moreover, by applying the uniform convergence of the
$\mu$-functions at limit stages of the partition construction, we
see that $\pi(t)_\beta = \mu(f_{(\beta,V^t_\beta)},t) = \lim_{\alpha
< \beta} \mu(f_{(\alpha,V^t_\alpha)},t) = \lim_{\alpha < \beta}
\pi(t)_\alpha$. This is enough to prove that $\pi(t) \in Z$.

We observe our first property of $\pi$, namely that it is
increasing. Let $s,t \in \Upsilon$ with $s \prec t$. We set $\gamma$
to be the least ordinal such that $\pi(s)_\gamma$ and
$\pi(t)_\gamma$ are not both defined and equal. If $\gamma = 0$
then, as $\mu$ is increasing, it follows that $\pi(s)_0 > \pi(t)_0$,
whence $\pi(s) < \pi(t)$. If $\gamma > 0$ then, by continuity,
$\gamma = \beta + 1$ for some $\beta$. By transfinite induction,
$V^s_\alpha = V^s_\alpha$ for all $\alpha \leq \beta$. Indeed,
$\mu(s) = -\pi(s)_0 = -\pi(t)_0 = \mu(t)$, so $V^s_0 = V^t_0$. If
$V^s_\alpha = U = V^t_\alpha$ and $\alpha < \beta$, set
$\lambda^s_\alpha = \lambda = \lambda^t_\alpha$. Remembering
property (2) of Lemma \ref{lambda}, we have
\begin{equation}
\label{lambdaeqn} \textstyle{\frac{1}{2}}\varepsilon\lambda(s) \;=\;
\pi(s)_{\alpha+1} - \pi(s)_\alpha \;=\; \pi(t)_{\alpha+1} -
\pi(t)_\alpha \;=\; \textstyle{\frac{1}{2}}\varepsilon\lambda(t)
\end{equation}
whence $\lambda(s) = \lambda(t)$ and $V^s_{\alpha+1} =
V^t_{\alpha+1}$. Limit stages of the induction follow by taking
intersections.

Now let $V^s_\beta = V = V^t_\beta$, $\lambda^s_\beta = \lambda
=\lambda^t_\beta$ and observe that $0_V \preccurlyeq s \prec t$.
There are two cases to consider: either $\pi(t)_{\beta+1}$ is
defined or it is not. First of all, we suppose that
$\pi(t)_{\beta+1}$ is defined and prove that $\pi(s) < \pi(t)$ in
this case. Indeed, if $\pi(s)_{\beta+1}$ is not defined then we are
done, as $\pi(t)$ strictly extends $\pi(s)$. On the other hand, if
$\pi(s)_{\beta+1}$ is defined then since $\pi(s)_{\beta+1} \neq
\pi(t)_{\beta+1}$ and $\lambda$ is decreasing, it must be that
$\pi(s)_{\beta+1} > \pi(t)_{\beta+1}$. Therefore $\pi(s) < \pi(t)$.

The other option is that $\pi(t)_{\beta+1}$ is undefined. In this
case, since $0_V \prec t$, it must be that $\lambda^t_\alpha(t) = 0$
for some $\alpha+1 < \beta+1$, by the definition of $\pi$. As
$\pi(t)_\beta$ is defined then, again by the definition of $\pi$, it
follows that $\alpha+1 = \beta$. Let $V^s_\alpha = U = V^t_\alpha$
and $\lambda^s_\alpha = \lambda^\prime = \lambda^t_\alpha$. Then by
Eqn.\ \ref{lambdaeqn} above, we have $\lambda^\prime(s) =
\lambda^\prime(t) = 0$, meaning $\pi(s)_{\beta+1}$ is not defined
either. Consequently, $\pi(s) = \pi(t)$.

We have established that $\pi$ is an increasing function. Now we
show that it is not constant on any ever-branching subset and, given
$t \in \Upsilon$, there are only finitely many $u \in t^+$ such that
$\pi(u) = \pi(t)$. To prove this claim, consider $t \in \Upsilon$
and the plateau $W = \setcomp{u \in [t,\infty)}{\pi(u) = \pi(t)}$.
If $W$ is the singleton $\{t\}$ then there is nothing to prove, so
we suppose that there exists some $u \in W$ with $t \prec u$. Let
both $\pi(t)$ and $\pi(u)$ be defined on $[0,\beta]$ and fix $V =
V^t_\beta$. In just the same way as above, we have that $V^t_\alpha
= V^u_\alpha$ whenever $\alpha \leq \beta$ and, in particular,
$V^u_\beta = V$. Observe that, as a consequence, $W \subseteq V$.
Moreover, just as above, as $\pi(u)_{\beta+1}$ is undefined and
$0_{V^u_\beta} \preccurlyeq t \prec u$, we have $\beta = \alpha+1$
for some $\alpha$. It follows that if we set $V^t_\alpha = U =
V^u_\alpha$ and $\lambda^t_\alpha = \lambda^\prime =
\lambda^u_\alpha$, then $\lambda^\prime(t) = \lambda^\prime(u) = 0$.
Now we can appeal to parts (2) and (3) of Lemma \ref{lambdavanish}
applied to $U$, $f_{(\alpha,U)}$, $\mu(f_{(\alpha,U)},\cdot)$ and
$\lambda^\prime$ to conclude that $V$ is finitely-branching and
contains no ever-branching subsets. As $W \subseteq V$, we are done.

We finish our appraisal of $\pi$ by showing that it does not admit
certain types of $Z$-bad points. First of all, if $\pi(t) \in Y$
then $t$ cannot be $Z$-bad for $\pi$. Indeed, by Lemma
\ref{convergenceinZ} part (1) and the fact that the elements of
$\ran \pi$ are uniformly bounded sequences, the only way that $t$
can be $Z$-bad for $\pi$ is if there are infinitely many $u \in t^+$
such that $\pi(u) = \pi(t)$. Now suppose that $\pi(t) =
(\pi(t)_\xi)_{\xi \leq \beta + 1} \in Z\backslash Y$, where $\beta$
is a limit ordinal. We prove that $t$ is not $Z$-bad for $\pi$. We
know already that $\pi(u) = \pi(t)$ for only finitely many $u \in
t^+$ so, for a contradiction, we must suppose that there is a
sequence of distinct points $(u_n) \subseteq t^+$ such that $\pi(t)
< \pi(u_n)$ and $\pi(u_n) \rightarrow \pi(t)$. We have that
$\pi(t)_\beta = \pi(t)_{\beta+1}$. Let $V = V^t_\beta$, where
$V^t_\beta$ is the unique element $V \in \mathscr{P}_\beta$
containing $t$, and let $f = f_{(\beta,V)}$. Observe that if
$\lambda$ is the function from Lemma \ref{lambda} associated with
$f$ and $V$ then, necessarily, $\lambda(t) = 0$. Indeed, by the
definition of $\pi$, we have $\frac{1}{2}\varepsilon \lambda(t) =
\pi(t)_{\beta+1} - \pi(t)_\beta$. By Lemma \ref{convergenceinZ} part
(3), there exist ordinals $\alpha_n < \beta$ such that $\alpha_n
\rightarrow \beta$, $\pi(u_n)_\xi = \pi(t)_\xi$ whenever $\xi <
\alpha_n$ and $\pi(u_n)_{\alpha_n} < \pi(t)_{\alpha_n}$. By
continuity and transfinite induction, $\alpha_n = \xi_n + 1$ for
some ordinals $\xi_n$ and $V^t_{\xi_n} = V^{u_n}_{\xi_n}$. Set $V_n
= V^t_{\xi_n}$ and $f_n = f_{(\xi_n,V_n)}$. As $\alpha_n \rightarrow
\beta$, it follows that $V = \bigcap_n V_n$ and the functions $f_n +
f_n(0_{V_n})\ind{(0_{V_n},t]}$ converge in norm to $f +
f(0_V)\ind{(0_V,t]}$. Moreover $\mu(f_n,u_n) = \pi(u_n)_{\xi_n} =
\pi(t)_{\xi_n} \rightarrow \pi(t)_\beta = \mu(f,t)$. Now choose
$\varphi_n \in \Czerok{\Upsilon}$ to satisfy $\supp \varphi_n
\subseteq (u_n, \infty)$ and $\norm{f_n +
f_n(0_{V_n})\ind{(0_{V_n},u_n]} + \varphi_n} \leq \mu(f_n,u_n) +
2^{-n} = \mu(f_n,t) + 2^{-n}$. As the $u_n$ are distinct, it follows
that $(f_n + f_n(0_{V_n})\ind{(0_{V_n},u_n]} + \varphi_n)$ converges
to $f + f(0_V)\ind{(0_V,t]}$ in the pointwise topology of
$\Czerok{\Upsilon}$. As $\Upsilon$ is scattered and this sequence is
norm-bounded, it converges in the weak topology too. Therefore
$\norm{f + f(0_V)\ind{(0_V,t]}} = \mu(f,t)$. However, by part (4) of
Lemma \ref{lambda}, the attainment of the infimum forces $\lambda(t)
> 0$, which is not the case. It follows that $t$ cannot be a $Z$-bad
point for $\pi$.

One case remains untreated. If $\pi(t) = (\pi(t)_\xi)_{\xi \leq
\beta + 1} \in Z\backslash Y$ and $\beta$ is not a limit ordinal, it
is possible that $t$ is $Z$-bad for $\pi$. Fortunately, by making an
adjustment to $\pi$ akin to that given after Lemma
\ref{convergenceinZ}, we can remove $Z$-bad points of this kind.
Given $x = (x_\xi)_{\xi \leq \beta} \in Z$, define
\[
\Phi(x)_\xi \;=\; \left\{
\begin{array} {l@{}l}
2x_0 & \quad\mbox{if } \xi = 0 \\
x_\xi + x_{\xi - 1} + 1 & \quad\mbox{if } \xi \mbox{ is a successor ordinal}\\
2x_\xi + 1 & \quad\mbox{otherwise}
\end{array} \right.
\]
for $\xi \leq \beta$. It is easy to establish that $\Phi$ takes
values in $Z$ and is strictly increasing. Set $\rho = \Phi \circ
\pi$. As $\Phi$ is strictly increasing, $\rho$ is increasing and, if
we consider Proposition \ref{plateaupartition}, partitions
$\Upsilon$ in exactly the same way as $\pi$. In particular, $\rho$
is not constant on any ever-branching subset of $\Upsilon$.  Again,
as $\Phi$ is strictly increasing, if $t$ is $Z$-bad for $\rho$ then
it is also $Z$-bad for $\pi$. Therefore, to prove that $\rho$ has no
$Z$-bad points, we suppose that $\pi(t) = (\pi(t)_\xi)_{\xi \leq
\beta + 1} \in Z\backslash Y$ and $\beta$ is not a limit ordinal. We
have that $\pi(t)_\beta = \pi(t)_{\beta+1}$ so, by the construction
of $\pi$, there exists an ordinal $\alpha$ such that $\beta = \alpha
+ 1$. Therefore, $\pi(t)_\alpha < \pi(t)_\beta$ and thus
$\rho(t)_\beta < \rho(t)_{\beta+1}$, giving $\rho(t) \in Y$. Again
by appealing to Lemma \ref{convergenceinZ} part (1), if $t$ is
$Z$-bad for $\rho$ then $\rho(u) = \rho(t)$ for infinitely many $u
\in t^+$. However, that would force $\pi(u) = \pi(t)$ for infinitely
many $u \in t^+$, and we have already established that this is
impossible.
\end{proof}

\begin{proof}[of Corollary \ref{latticeggivesr}]
If $\normdot$ is a lattice norm then, by part (1) of Lemma
\ref{infattaining}, the infima in the definition of the
$\mu$-functions are always attained. It follows that the
$\lambda$-functions of Lemma \ref{lambda} never vanish. Now, we
prove that in this case, the map $\pi$ defined in the proof of
Theorem \ref{newgnec} is $Y$-valued and strictly increasing. Indeed,
if we return to the point where we prove that $\pi(t) \in Z$, we see
that, as the $\lambda$-functions never vanish, $\pi(t)_\alpha <
\pi(t)_{\alpha+1}$ whenever $\alpha+1 \leq \beta$. Consequently
$\pi(t) \in Y$. To show that $\pi$ is strictly increasing, we let $s
\prec t$ and return to the point in the proof where $\pi$ is shown
to be increasing, specifically, where $\gamma$ is defined. If
$\gamma = 0$ then we are done. Otherwise, $\gamma = \beta+1$ for
some $\beta$. Since the $\lambda$-functions never vanish, it is
impossible that $\pi(t)_{\beta+1}$ is undefined, therefore $\pi(s) <
\pi(t)$. This proves that $\Upsilon \preccurlyeq Y$. The second
statement of Corollary \ref{latticeggivesr} holds because the
strictly convex dual norm constructed in Theorem \ref{dualrotundthm}
is a lattice norm.
\end{proof}

We finish the section with a proof of Proposition \ref{haydoncor}.
It will help to introduce a useful game-theoretic characterisation
of Baire trees \cite{haydon:95}. Players \textbf{A} and \textbf{B}
take turns to nominate elements of a tree $\Upsilon$, beginning with
$t_0$ played by \textbf{B}. In general, \textbf{A} follows $t_{2n}$
with $t_{2n+1} \succcurlyeq t_{2n}$, and \textbf{B} responds with
$t_{2n+2} \succcurlyeq t_{2n+1}$. The game is won by \textbf{B} if
the sequence $(t_n)$ has no upper bound in $\Upsilon$. The tree
$\Upsilon$ is Baire if and only if \textbf{B} has no winning
strategy in this so-called $\Upsilon$\textit{-game}. Using this
game, it is possible to prove the following result.

\begin{prop}[({Haydon \cite[Proposition 1.4]{haydon:95}})]
\label{realsubmissive} If $\Upsilon$ is Baire and
$\mapping{\rho}{\Upsilon}{\real}$ is increasing, then there exists
$t \in \Upsilon$ such that $\rho$ is constant on the wedge
$[t,\infty)$.
\end{prop}

One trivial consequence of Proposition \ref{realsubmissive} is that
if the increasing map $\mapping{\rho}{\Upsilon}{\real}$ is not
constant on any ever-branching subset then $\Upsilon$ contains no
ever-branching Baire subsets. Indeed, if $E \subseteq \Upsilon$ were
ever-branching and Baire then, by Proposition \ref{realsubmissive},
we could find $t \in E$ such that $\rho$ is constant on $[t,\infty)
\cap E$, which is an ever-branching subset of $\Upsilon$. We observe
that the same holds if we replace $\real$ with any linear order
$\Sigma$ satisfying the statement of Proposition
\ref{realsubmissive}. Therefore, to establish Proposition
\ref{haydoncor}, it is enough to prove the following result.

\begin{prop}
\label{zsubmissive} If $\Upsilon$ is Baire and
$\mapping{\rho}{\Upsilon}{Z}$ is increasing, then there exists $t
\in \Upsilon$ such that $\rho$ is constant on $[t,\infty)$.
\end{prop}

\begin{proof}
The following order will be used in this and a subsequent proof.
Define
\[
Z_0 \;=\; \setcomp{x = (x_\alpha)_{\alpha \leq \beta} \in Z}{x
\subseteq [0,1]\mbox{, }x_0 = 0\mbox{ and }\beta\mbox{ is a limit
whenever }x_\beta = 1}.
\]
By considering the map $\Theta$, introduced after Lemma
\ref{convergenceinZ}, we observe that $Z \preccurlyeq Z_0$ and,
accordingly, we can assume that our increasing function $\rho$ takes
values in $Z_0$.

We show that $\rho$ is constant on some wedge of $\Upsilon$ by
playing the $\Upsilon$-game with a particular strategy for
\textbf{B}. Given $u \in \Upsilon$ and an ordinal $\alpha$, we call
$(\alpha,u)$ a \textit{fixed pair} if $\rho(v)_\xi$ is defined and
equal to $\rho(u)_\xi$ whenever $v \in [u,\infty)$ and $\xi \leq
\alpha$. If $(\alpha,u)$ is fixed, $v \in [u,\infty)$ and $\xi \leq
\alpha$, then $(\xi,v)$ is also fixed. Let \textbf{B} play arbitrary
$t_0$ as the first move and put $\alpha_0 = 0$. Note that $(0,t_0)$
is fixed. Now suppose that $n \geq 1$ and that moves $t_0
\preccurlyeq t_1 \preccurlyeq \ldots \preccurlyeq t_{2n-1}$ have
been played alternately by \textbf{B} and \textbf{A}. We choose the
next move $t_{2n}$ played by \textbf{B}, together with $\alpha_n$,
in the following manner. Let
\[
r_n \;=\; \sup \setcomp{\rho(u)_\alpha}{u \succcurlyeq
t_{2n-1}\mbox{ and }(\alpha,u)\mbox{ is a fixed pair}}.
\]
Let \textbf{B} choose fixed $(\alpha_n,t_{2n})$ such that $t_{2n}
\succcurlyeq t_{2n-1}$ and $\rho(t_{2n})_{\alpha_n}
> r_n - 2^{-n}$. This strategy does not guarantee a win for \textbf{B}, so
there exist moves $(t_{2n+1})$ of \textbf{A} such that $(t_n)$ has
an upper bound $u \in \Upsilon$. If $\alpha = \sup \alpha_n$, we see
that $(\alpha,u)$ is fixed. This follows by continuity and the fact
that $(\alpha_n,u)$ is fixed for all $n$.

If $\rho(v)_{\alpha+1}$ is not defined for any $v \succcurlyeq u$
then $\rho$ takes constant value $\rho(u)$ on $[u,\infty)$, and we
are done. Suppose instead that $\rho(v)_{\alpha+1}$ exists for some
$v \succcurlyeq u$. Because $(\alpha,v)$ is fixed and $\rho$ is
increasing, the real-valued map $\rho(\cdot)_{\alpha+1}$ must be
decreasing on $[v,\infty)$. As the forcing-open set $[v,\infty)$ is
Baire, by Proposition \ref{realsubmissive}, there exists $w
\succcurlyeq v$ such that $\rho(\cdot)_{\alpha+1}$ is constant on
$[w,\infty)$, and it follows that $(\alpha+1,w)$ is a fixed pair. We
note that the inequalities
\[
r_n - 2^{-n} \;<\; \rho(t_{2n})_{\alpha_n} \;=\; \rho(w)_{\alpha_n}
\;\leq\; \rho(w)_\alpha \;\leq\; \rho(w)_{\alpha+1} \;\leq\; r_n
\]
hold for all $n$, and conclude that $\rho(w)_{\alpha+1} =
\rho(w)_\alpha$. Consequently, by the definition of elements of $Z$,
$\rho$ takes constant value $\rho(w)$ on $[w,\infty)$.
\end{proof}

\section{Examples}
\label{examples}

In this section, we prove Proposition \ref{sigmaYnoG} and present
Example \ref{Ggap}. Before giving the proof of Proposition
\ref{sigmaYnoG}, we make an observation about embeddability and
$Z$-bad points that is analogous to Proposition \ref{ratbadpoints}.

Given a tree $\Upsilon$, let $\Upsilon \preccurlyeq Z$ and suppose
that there is an increasing function $\mapping{\rho}{\Upsilon}{Z}$
with no $Z$-bad points. We claim that if this is the case then
$\Upsilon \preccurlyeq Y$. In order to prove this claim, we
introduce the following algebraic operation on $Z$. Recall the order
isomorphism $\mapping{\theta}{\real}{(0,1)}$, fixed after Lemma
\ref{convergenceinZ}. For $x = (x_\xi)_{\xi \leq \alpha}$ and $y =
(y_\xi)_{\xi \leq \beta}$ of $Z$, define $x\cdot y$ for $\xi \leq
\max\{\alpha,\beta\}$ by
\[
(x\cdot y)_\xi \;=\; \left\{
\begin{array} {l@{}l}
\theta^{-1}(\theta(x_\xi)\theta(y_\xi)) & \quad\mbox{if } \xi \leq \min\{\alpha,\beta\} \\
x_\xi & \quad\mbox{if } \alpha < \xi \leq \beta\\
y_\xi & \quad\mbox{if } \beta < \xi \leq \alpha
\end{array} \right.
\]
where $\theta(x_\xi)\theta(y_\xi)$ is an ordinary real product. We
leave the reader with the simple task of verifying that $\cdot$ is a
semigroup operation on $Z$ that respects the order; in other words,
if $x \leq y$ and $u \leq v$ then $x\cdot u \leq y\cdot v$ and,
moreover, the third inequality is strict if either of the first two
are. Now, let the increasing function $\mapping{\nu}{\Upsilon}{Z}$
have no $Z$-bad points and suppose $\mapping{\tau}{\Upsilon}{Z}$ is
strictly increasing. As $\cdot$ respects order, it follows that the
pointwise product $\pi = \nu\cdot\tau$ is strictly increasing and
has no $Z$-bad points. By Lemma \ref{convergenceinZ}, any element of
$Z$ can be approached from above by a strictly decreasing sequence.
Therefore, as $t \in \Upsilon$ is not a $Z$-bad point for $\pi$,
there exists $\pi^*(t) \in Z$ such that $\pi(t) < \pi^*(t) \leq
\pi(u)$ whenever $u \in t^+$.  Finally, since $Y$ is dense in $Z$,
we can pick $\rho(t) \in Y$ between $\pi(t)$ and $\pi^*(t)$; the
resulting function $\rho$ is strictly increasing.

\begin{proof}[of Proposition \ref{sigmaYnoG}]
In the light of Theorem \ref{sigmanoembed} and our observation
above, all we need to do is prove that $\sigma Y \preccurlyeq Z$.
Recall the order $Z_0$ from the proof of Proposition
\ref{zsubmissive}. As $Z \preccurlyeq Z_0$, elements of $\sigma Y$
can and are considered as subsets of $Z_0$. Our proof that $\sigma Y
\preccurlyeq Z$ rests on the claim that $Z_0$ is Dedekind complete;
that is, each subset of $A$ of $Z_0$ has a least upper bound,
denoted by $\sup A$.

For now, we assume that this claim holds and define a strictly
increasing map $\mapping{\rho}{\sigma Y}{Z}$. Given $A \in \sigma
Y$, treated as a subset of $Z_0$, let $\rho(A) = \sup A$ if $\sup A
\in Z_0\backslash Y$ or if $A$ has no greatest element, and let
$\rho(A) = (\sup A,2)$ otherwise. Here, $(x,2)$ denotes the sequence
obtained by extending $x \in Z_0 \cap Y$ by a single element, namely
$2$. Observe that if $x \in Z_0 \cap Y$, $y \in Z_0$ and $x < y$
then $(x,2) < y$ because every element of $y$ is strictly less than
$2$. Let $A,B \in \sigma Y$ satisfy $A \prec B$. If $\sup A < \sup
B$ then $\rho(A) < \sup B \leq \rho(B)$. Alternatively, if $\sup A =
\sup B$ then $B = A \cup \{\sup A\}$; indeed, if $x \in B\backslash
A$ then $\sup A \leq x \leq \sup B = \sup A$. In particular, $B$ has
greatest element $\sup A \in Y$, whereas $A$ has no greatest
element. Therefore $\rho(A) = \sup A < (\sup A,2) = \rho(B)$. This
proves that $\rho$ is strictly increasing.

To finish, we define $\sup A$ for $A \subseteq Z_0$. If $A$ is empty
then its least upper bound is the one-element sequence $(0)$. From
now on, we assume that $A$ is non-empty and has no greatest element.
Taking our cue from the proof of Proposition \ref{zsubmissive},
given an ordinal $\alpha$ and $x \in A$, we will call $(\alpha,x)$ a
\textit{fixed pair} if $x_\xi$ and $y_\xi$ are both defined and
equal whenever $y \in A$, $x \leq y$ and $\xi \leq \alpha$. If
$(\alpha,x)$ is fixed, $y \in A$, $x \leq y$ and $\xi \leq \alpha$,
then $(\xi,y)$ is also fixed. Now let $\beta$ be minimal, subject to
the condition that there is no fixed pair $(\beta,x)$. As A is
non-empty and $(0,x)$ is fixed whenever $x \in A$, it follows that
$\beta > 0$. We define a sequence $z = (z_\alpha)_{\alpha \leq
\beta}$. If $\alpha < \beta$, let $z_\alpha = x_\alpha$, where
$(\alpha,x)$ is some fixed pair. By the nature of fixed pairs, this
is well-defined. If $\beta$ is a limit, let $z_\beta = \sup_{\alpha
< \beta} z_\alpha$. Instead, if $\beta = \alpha + 1$ for some
$\alpha$ then, as $A$ has no greatest element, there exists a fixed
pair $(\alpha,x)$, such that $x_\beta$ is defined. Let $z_\beta$ be
the infimum of all such $x_\beta$. It is easy to verify that $z \in
Z_0$; it can be that $z_\beta = 1$, but only if $\beta$ is a limit
ordinal. We omit the pedestrian task of proving that $z$ is the
least upper bound of $A$.
\end{proof}

Our last task is to show that there is a tree $\Psi$ satisfying the
condition of Theorem \ref{newgnec} but not that of Theorem
\ref{gsuff}. Before doing so, we must make some remarks. Recall the
plateau partitions of Proposition \ref{plateaupartition} and note
the following slightly reworded version of a result from
\cite{smith:05c}.

\begin{prop}[({Smith \cite[Corollary 3]{smith:05c}})]
\label{noebcor} Suppose that $\Upsilon$ is a tree, $\Sigma$ a linear
order, and $\mapping{\rho}{\Upsilon}{\Sigma}$ an increasing function
that is not constant on any ever-branching subset of $\Upsilon$.
Then there exists an increasing function
$\mapping{\pi}{\Upsilon}{\Sigma \times \omega}$, such that the
plateau partition $\mathscr{P}$ of $\Upsilon$ with respect to $\pi$
consists solely of linearly ordered subsets.
\end{prop}

Let $\Upsilon$, $\Sigma$, $\pi$ and $\mathscr{P}$ be as in
Proposition \ref{noebcor} and, moreover, let us suppose that
$\Upsilon$ admits no uncountable linearly ordered subsets. In this
case, each $V \in \mathscr{P}$ identifies with a finite or countable
ordinal and, therefore, there exists a strictly increasing function
$\mapping{\pi_V}{V}{\rat}$. It is apparent that the function
$\mapping{\tau}{\Upsilon}{\Sigma \times \omega \times \rat}$,
defined by $\tau(t) = (\pi(t),\pi_{V_t}(t))$, where $V_t$ is the
unique element of $\mathscr{P}$ containing $t$, is strictly
increasing. As $\omega \times \rat \preccurlyeq \rat$, it follows
that $\Upsilon \preccurlyeq \Sigma \times \rat$.

\begin{example}
\label{Ggap} Observe that $Y$ has cardinality continuum
$\continuum$. If $A \in \sigma Y$ then $A^+$ identifies with the set
$u(A)$ of all upper bounds of $A$ and, thus, has cardinality
$\continuum$ if $u(A)$ is non-empty. Fix a well-order $\sqsubseteq$
of $Y$, and let $\Psi = \sigma Y \times \continuum$. We order $\Psi$
by declaring that $(A,\alpha) \preccurlyeq (B,\beta)$ if and only if
either $A = B$ and $\alpha \leq \beta$, or if $A \prec B$ and
$\alpha$ is no greater than the order type of $\setcomp{x \in
u(A)}{x \sqsubset \min (B\backslash A,\leq)}$, with respect to
$\sqsubset$.

With respect to this order, each element of $\Psi$ has between one
and two immediate successors. Indeed, if $(A,\alpha) \in \Psi$ then
$(A,\alpha + 1)$ is always an immediate successor. If $u(A)$ is
non-empty then $(A \cup \{y\},0)$ is also such a successor, where $y
\in u(A)$ and $\setcomp{x \in u(A)}{x \sqsubset y}$ has order type
$\alpha$. The set $\sigma Y \times \{0\}$ is a natural copy of
$\sigma Y$ inside $\Psi$ that is closed with respect to the interval
topology.

Now, by Proposition \ref{sigmaYnoG}, there exists a strictly
increasing map $\mapping{\pi}{\sigma Y}{Z}$. Define
$\mapping{\rho}{\Psi}{Z}$ by $\rho(A,\alpha) = \pi(A)$. By
Proposition \ref{plateaupartition}, the plateau partition of $\Psi$
with respect to $\rho$ consists exactly of the sets
$\setcomp{(A,\alpha)}{\alpha < \continuum}$, where $A \in \sigma Y$.
Therefore, $\rho$ is not constant on any ever-branching subset.
Because the number of immediate successors of any element of $\Psi$
is at most two, $\rho$ has no $Z$-bad points either. Therefore
$\Psi$ satisfies the condition of Proposition \ref{newgnec}.

On the other hand, there exists no increasing $Y$-valued function on
$\Psi$ that is not constant on any ever-branching subset. Indeed, if
there were such a function, by considering its restriction to
$\sigma Y \times \{0\}$, there would be a map $\mapping{\tau}{\sigma
Y}{Y}$, also not constant on any ever-branching subset. However, by
following a similar argument to that given after Proposition
\ref{realsubmissive}, being $Z$-embeddable, $\sigma Y$ has no
perfect Baire subsets. In particular, $\sigma Y$ does not contain a
copy of $\wone$. Therefore, by Proposition \ref{ybetaembed} and the
remarks following Proposition \ref{noebcor}, we would have $\sigma Y
\preccurlyeq Y \times \rat \preccurlyeq Y$ which, by Theorem
\ref{sigmanoembed}, is impossible.
\end{example}

We recall Problem \ref{gateauxconj} and conjecture that
$\Czerok{\Psi}$ admits a G\^{a}teaux norm. The G\^{a}teaux norms
presented in \cite{smith:05c} are built by combining norms obtained
from existing techniques, namely the Fr\'{e}chet norms of Talagrand
and Haydon, and norms with strictly convex duals. In the author's
opinion, if Problem \ref{gateauxconj} is to be resolved positively,
we require a method of constructing G\^{a}teaux norms on $\Ck{K}$
spaces that unifies these techniques on a more fundamental level.

\bibliographystyle{amsplain}

\end{document}